\documentstyle[12pt]{article}
\begin{document}
\title{Mixed means of commutators \\of  central integral means and $CMO$   }
\author{{ Shunchao Long and Jian Wang }    }
 \date{}
\maketitle
\begin{center}
\begin{minipage}{120mm}
\vskip 0.2in {{\bf }~~ In this paper we obtain some mixed means and
weighted $L^p$ estimates for the commutators generating $r$ order
central integral means operators  with
  $CMO$ functions.
 }
\end{minipage}
\end{center}
\vskip 0.2in \baselineskip 0.2in
\par \section*{{ \bf  1.INTRODUCTION}
} ~~~~  \par Let $T$ be a sublinear operator and $b(x)$ a function,
the  commutator $[T, b]$  is defined by
$$[T, b]f(x) = T(fb)(x)-b(x)Tf(x).$$

A famous result of Coifman,  Rochberg and Weiss stated  that the
commutators of Calderon-Zygmund singular integral operators and
$BMO$ functions   are bounded on $L^p$ for $1<p<\infty$ ( see [5]).
Since then, the $L^p$ estimates and applications of
 these type commutators were studied by many authors (see [2,21,19,4] and [3,15,20,23]).


\par Recently the authors of this paper proved  the $L^p$-boundedness of the commutators of Hardy operators
and  one-side $CMO$ functions  for $1<p<\infty$ in [18]. And this
result was extended in [16,11].
 In this paper we  establish  some mixed means estimates and weighted
 $L^p$-estimates of the commutators of central integral means and $ CMO$ functions, and this extends the boundedness results in [18,11].

  \par Let $1\leq p<{\infty}$. Denote
\[
\sup\limits_{B\subset {\bf R}^n} \left( \frac {1}{|B|}\int_{B}|b(x)-b_{B}|^pdx \right)^{1/p}
=\left\{
\begin{array}{cc}
\|b\|_{  BMO ^p}, & {\rm if} ~ B ~{\rm are~ arbitrary~ balls},\\
\|b\|_{  CMO ^p}, & {\rm if} ~ B ~{\rm are~  balls~ centraled~ at~
origin},
\end{array}
\right.
\]
where $b_{B}=\frac {1}{|B|}\int_{B}b(x)dx $. If $\|b\|_{  BMO
^p}<\infty,$ we say $b\in   BMO  ^p,$ the Bounded Mean Oscillation.
 It is well-known that $  BMO ^p=  BMO ^1=  BMO  $ for all $1\leq p<{\infty}$. If $\|b\|_{  CMO  ^p}<\infty,
 $ we say $b\in   CMO  ^p,$ the Central Mean Oscillation. Obviously, if $1\leq p<q<{\infty},$ then
  $  CMO  ^q\subset   CMO  ^p$, and $\|b\|_{  CMO  ^p}\leq \|b\|_{  CMO  ^q}.$
\par We say $b\in  CMO$, if
$ \|b\|_{  CMO   }=\sup _{1\leq p< \infty  }\|b\|_{ CMO  ^p}<\infty$.
\par The spaces $ CMO ^p   $   bear a simple relationship with $BMO$: $g \in  BMO  $ precisely
when $g$ and all of its translates belong to $CMO ^p$ uniformly a.e., so
 $${\rm the~ classical~} BMO {\rm~space} ~\subset   CMO  ^p~{\rm for~ all ~}1\leq p<{\infty}.$$
Many precise analogies exist between $CMO^p$ and $BMO$ from the
point of view of real Hardy spaces, for example, the duality:
$CMO^{p'}, p'=p/(p-1),$ is the dual of the Beurling-Herz-Hardy
spaces $HA^p, 1 \leq p< \infty $, that is then analogous to the $H^1
\leftrightarrow BMO$, (see [ 11-12 ] );
 the constructive decomposition (see [ 17 ]);
   and so on.
 \par If $R>0,$ denote $B(R)=B(0,R)$ be the ball in ${\bf R}^n$ centered at origin and of radius $R$.
Let $r, \alpha \in {\bf R}, f\in L^r_{loc}{({\bf R}^n, |x|^{n(\alpha -1)})}$,    the central integral mean
  of order $r$, with the power weight, of $f$, be defined  by
\begin{eqnarray*}
 M_{r }(f,  \alpha )(|y|)=\left[\frac{1}{|B(|y|)|^{\alpha}}
\int_{B(|y|)}|B(|x|)|^{\alpha-1} |f(x)|^rdx\right]^{1/r},
 \end{eqnarray*}
and the companion mean   of order $r$, with the power weight, of $f$, by
\begin{eqnarray*}
 M^*_{r} (f,  \alpha  )(|y|)=\left[\frac{1}{|B(|y|)|^{\alpha}}
\int_{{\bf R}^n\setminus B(|y|)}|B(|x|)|^{\alpha-1} |f(x)|^rdx\right]^{1/r}.
 \end{eqnarray*}
\par The properties and applications of these types of integral means can be found in many
literatures. Firstly,
  the limits of   $(M_{2 }(f,  1)(y) )^{2}=(1/2y)\int^y_{-y}|f|^2(y>0$,
  the one-dimensional case),  were used to study the almost periodic functions,
   and the spectrum and ergodicity of sample paths of certain stochastic processes
    by Wiener  in [ 22 ].
  Secondly, the functions spaces of bounded integral mean of $r$ order, introduced firstly by Beurling,
$$B^r= B^{r, \infty} =\{f: M_{r }(f,  1)(|y|) \in L^{\infty} \},$$
   (both homogeneous spaces ($|y|>0$) and non-homogeneous spaces $(|y|>1))$, together
   with their corresponding Beurling algebra $A^r$ and the Hardy space $HA^r$
   [ cf, 1, 12, 6  ]
    had  rich contents;
Thirdly, the Hardy type inequalities   associating with $M_{1 }(f,
1)(|y|)$ and $M^*_{1 }(f,  1 )(|y|)$ generalized the classical ones
to $n$-dimension ball [ 9-10, 7  ]. And the mixed means inequalities
were used to derive the generalizing Hardy and Levin-Cochran-Lee
type inequalities in [ 7 ] (see also [8]).
\par We state the mixed means inequalities as following:
    \par {\bf Theorem 1 } ([ 7: Theorem 5]).~~ Let $r ,s, R, \alpha , \gamma \in {\bf
 R},$
  and let $  r<s,  r,s\not= 0, R>0 (f\not= 0  $ in the case of $r<0$). Then,
\begin{eqnarray}
  M_{s}((M_{r  }(f,  \alpha ),\gamma)(R)\leq M_{r}((M_{s  }(f,  \gamma ),\alpha)(R),
\end{eqnarray}
 \begin{eqnarray}
  M^*_{s}((M^*_{r  }(f,  \alpha ),\gamma)(R)\leq M^*_{r}((M^*_{s  }(f,  \gamma ),\alpha)(R).
\end{eqnarray}
 \par From Theorem 1, we can obtain the $L^p$-boundedness estimates of these
 types of
 integral means above as following.

\par {\bf Theorem 2.~~} Let $r ,s,   \alpha , \gamma \in {\bf
 R},$
  and let $  r<s,  r \not= 0,  s>0 (f\not= 0  $ in the case of $r<0$).   Then,
  \begin{eqnarray}
   \int_{{\bf R}^n}|B(|y|)|^{\gamma-1}(M_{r }(f,  \alpha
   )(|y|))^sdy
 \leq \frac {1}{(\alpha-\gamma r/s)^{s/r}}  \int_{{\bf R}^n} |B(|y|)|^{\gamma-1}|f
 (y)|^sdy
\end{eqnarray}
if $\alpha-\gamma r/s>0,$ and
\begin{eqnarray}
   \int_{{\bf R}^n}|B(|y|)|^{\gamma-1}(M^*_{r }(f,  \alpha
   )(|y|))^sdy
 \leq \frac {1}{( \gamma r/s-\alpha)^{s/r}}   \int_{{\bf R}^n} |B(|y|)|^{\gamma-1}|f
 (y)|^sdy
\end{eqnarray}
if $\alpha-\gamma r/s<0.$
 \par {\bf   Proof  .}
  For $y\in B(R),$ $\int_{B(|y|)}|B(|x|)|^{\gamma-1}|f(x)|^s dx\leq
  \int_{B(R)}|B(|x|)|^{\gamma-1}|f (x)|^s dx$   and $\int_{B(R)}
  |B(|y|)|^{\alpha-1-\gamma r/s}dy
   =\frac {1}{\alpha-\gamma r/s} |B(R)|^{\alpha-\gamma r/s}
$ when $\alpha-\gamma r/s>0,$ and using (1), we have
\begin{eqnarray*}
& &  \int_{B(R)}|B(|y|)|^{\gamma-1}(M_{r  }(f,  \alpha )(|y|))^sdy
\\
& &
\leq {|B(R)|^{\gamma -s\alpha /r}} \left[ \int_{B(R)}
|B(|y|)|^{\alpha-1} \left(\frac{1}{|B(|y|)|^{\gamma}}
\int_{B(|y|)}|B(|x|)|^{\gamma-1}|f (x)|^s dx\right)^{r/s}
dy\right]^{s/r}
\\
& & \leq \frac {1}{(\alpha-\gamma r/s)^{s/r}}
\int_{B(R)}|B(|x|)|^{\gamma-1}|f (x)|^s dx,
\end{eqnarray*}
(3) follows by taking the lim$_{R\rightarrow \infty}$. And (4) is
the consequence of (2), derived by the same technique as (3) from
(1), except for taking the lim$_{R\rightarrow 0}$.

\par {\bf Remark} ~~It is easy to see that $L^\infty \subset B^{r,\infty}$
for $0<r<\infty$ (see also [6,12-13]). Let
$$  B^{r,s}(\alpha, \gamma) =\{f: M_{r }(f,  \alpha )(|y|) \in L^{s}_{|x|^{n\gamma}} \}.$$
Then   $r, s,   \alpha , \gamma  $ are under the conditions of
Theorem 2 and $\alpha-\gamma r/s>0$. Thus, we have
$$L^{s}_{|x|^{n\gamma}}\subset B^{r,s}(\alpha, \gamma).$$

\par Further, we can obtain the boundedness estimates of the commutators generating $r$ order central integral means operators  and
  $CMO$ functions.
\par From the point of view on   Hardy spaces, the central integral means  bear  some  relationships with $ CMO ^p$. In fact,
  $(A^p)^*=B^{p'}$ and $(HA^p)^*=CMO^{p'}$.

\par
\section * {\bf 2. COMMUTATOR THEOREMS   }
~~~~ Let $r>0, \alpha \in {\bf R}, f\in L^r_{loc}{({\bf R}^n,
|x|^{n(\alpha -1)})}$, and $b$ be local integral functions on ${\bf
R}^n  $. We define the integral mean commutators  of order $r$, with
the power weight, by
\begin{eqnarray*}
 M_{r,b }(f,  \alpha )(|y|)=\left[\frac{1}{|B(|y|)|^{\alpha}}
\int_{B(|y|)}|B(|x|)|^{\alpha-1}(|b(x)-b(y)| |f(x)|)^rdx\right]^{1/r},
 \end{eqnarray*}
and the companion mean commutators  of order $r$, with the power weight, by
\begin{eqnarray*}
 M^*_{r,b }(f,  \alpha  )(|y|)=\left[\frac{1}{|B(|y|)|^{\alpha}}
\int_{{\bf R}^n\setminus B(|y|)}|B(|x|)|^{\alpha-1}(|b(x)-b(y)| |f(x)|)^rdx\right]^{1/r}.
 \end{eqnarray*}
 When $r\ge 1$, by the Minkowski inequality, it is easy to calculate that,
  $$
|[M_r(\bullet,\alpha),b]f(|y|)|\stackrel {def}{=}
|M_r(bf,\alpha)(|y|)-b(|y|)M_r(f,\alpha)(|y|)|\leq M_{r,b }(f,
\alpha )(|y|),
$$
$$
|[M^*_r(\bullet,\alpha),b]f(|y|)|\stackrel {def}{=}
|M^*_r(bf,\alpha)(|y|)-b(|y|)M^*_r(f,\alpha)(|y|)|\leq M^*_{r,b
}(f,  \alpha )(|y|).
$$
\par  Mixed-means inequalities of commutators of integral means and $CMO$ functions:
\par {\bf Theorem 3.~~}Let $ r ,s, R, \alpha , \gamma \in {\bf R}, s>r>0,
R>0.$ Suppose that $
      b \in CMO $ and $ f  $ is local  bounded   functions.   Then,
\begin{eqnarray}
  M_{s}((M_{r,b  }(f,  \alpha ),\gamma)(R)\leq c_1 \|b\|_{  CMO   }M_{r}((M_{s  }(f,  \gamma ),1)(R)
\end{eqnarray}
 if $\alpha >1$; and
\begin{eqnarray}
  M^*_{s}((M^*_{r,b  }(f,  \alpha ),\gamma)(R)\leq c_1 \|b\|_{  CMO   }M^*_{r}((M^*_{s  }(f,  \gamma ),1)(R)
\end{eqnarray}
    if $\alpha <1,$ where $c_1 = [2^{n|\alpha |}2^{n|\alpha -1|}3^{r  }\times4\times 2^{2nr }\sum\limits_{k=0}^{ \infty}
  2^{-kn|\alpha-1| }k^{r }]^{1/r}$.
  \par Weighted $L^p-$estimates of commutators of integral means and $CMO$ functions:
\par {\bf Theorem 4.~~} Let $r ,s,   \alpha , \gamma \in {\bf R}, s>r>0 ;
     b \in CMO, f $ is local bounded functions,
   then,
\begin{eqnarray}
   \int_{{\bf R}^n}|B(|y|)|^{\gamma-1}(M_{r,b }(f,  \alpha )(|y|))^sdy
 \leq c_2 \|b\|_{  CMO   }^s \int_{{\bf R}^n} |B(|y|)|^{\gamma-1}|f (y)|^sdy  ,
\end{eqnarray}
if $\alpha >1 $ and $\gamma <s/r,$ and
\begin{eqnarray}
   \int_{{\bf R}^n}|B(|y|)|^{\gamma-1}(M^*_{r,b }(f,  \alpha )(|y|))^sdy  \leq c_2  \|b\|_{  CMO   }^s \int_{{\bf R}^n}
    |B(|y|)|^{\gamma-1}|f (y)|^sdy  ,
\end{eqnarray}
 if $\alpha <1$ and $\gamma >s/r,$ where $c_2 = \frac {1}{|1-\gamma r/s|^{ s/r}}
 [2^{n|\alpha |}2^{n|\alpha -1|}3^{r  }\times4\times 2^{2nr }\sum\limits_{k=0}^{ \infty}
  2^{-kn|\alpha-1| }k^{r }]^{s/r}$.

 \par {\bf Proof of Theorem 3}~~ Let us prove (5). Let $2^{N-1}<R<2^{N},$ denote
\begin{eqnarray}
B_i=\left \{\begin{array}{cc}
B(2^i), & {\rm if}~ i\leq N-1,\\
B(R), & {\rm if}~ i= N,
\end{array}
\right. , C_i = B_i \setminus B_{i-1}, i=-\infty, ... , N;
\end{eqnarray}
  if $x\in C_i$, denote
 \begin{eqnarray}
{\overline {B_j}}=\left \{
\begin{array}{cc}
B_j, & {\rm if}~ j\leq i-1,\\
B(|x|), & {\rm if}~ j= i,
\end{array}
\right.,
 \overline {C_j} = \overline {B_j} \setminus \overline
{B_{j-1}}, j=-\infty, ... , i.
\end{eqnarray}
  Then, we have
\[
h(x)=\left[  (M_{r,b  }(f,  \alpha ) (|x|)\right]^r
=\frac{1}{|B(|x|)|^{\alpha}}
\sum\limits_{j=-\infty}^{i}
 \int_{\overline {C_j}}|B(|y|)|^{\alpha -1}
(|b(x)-b(y)| |f(y)|)^rdy,
\]
and
\begin{eqnarray}
  \left[ M_{s}((M_{r,b  }(f,  \alpha ),\gamma)(R)\right]^s
 = \frac{1}{|B(R)|^{\gamma}} \sum\limits_{i=-\infty}^{N}
 \int_{C_i}|B(|x|)|^{\gamma -1}
 \left[ h(x)  \right]^{s/r}dx.
\end{eqnarray}
If $x\in   {C_i}$ or $\overline {C_i}$, using $|B(R)|=R^n|B(1)|,$ we
have
\begin{eqnarray}
|B(|x|)|^{\alpha}\leq
\left \{\begin{array}{cc}
|B(1)|^{ \alpha } 2^{in\alpha } , & {\rm if}~ \alpha\geq 0,\\
|B(1)|^{ \alpha } 2^{(i-1)n\alpha } , & {\rm if}~ \alpha< 0,
\end{array}
\right\}
\leq |B(1)|^{ \alpha } 2^{n|\alpha |}2^{in\alpha },
\end{eqnarray}
and for $r>0 $ ,
\begin{eqnarray}
|a+b+c|^{r}\leq
\left \{\begin{array}{cc}
 |a |^{r}+| b |^{r}+| c|^{r} , & {\rm if}~ 0<r\leq 1,\\
3^{r-1 } (|a |^{r}+| b |^{r}+| c|^{r}), & {\rm if}~ 1< r,
\end{array}
\right\} \leq 3^{r } (|a |^{r}+| b |^{r}+| c|^{r}).
\end{eqnarray}
The first inequality is obvious when $0<r\leq 1,$ from the property
of convex function [14] when $r>1
 $ since $  g(x)=x^r$ is convex function. Noticing that $|b(x)-b(y)|\leq |b(x)-b_{
{B_i}}|+|b(y)-b_{\overline {B_j}}|+|b_{  {B_i}}-b_{\overline
{B_j}}|$ , using (13) and (12), we have
 \begin{eqnarray*}
h(x)& \leq & 2^{n|\alpha |}2^{n|\alpha -1|}3^{r }
  \frac{1}{2^{in\alpha }|B(1)|}\sum\limits_{j=-\infty}^{i}
  2^{jn(\alpha-1) }
 \int_{\overline {C_j}} |b(x)-b_{  {B_i}}|^{r }|f(y)|^rdy
\\
& & +2^{n|\alpha |}2^{n|\alpha -1|}3^{r  }
 \frac{1}{2^{in\alpha }|B(1)|}\sum\limits_{j=-\infty}^{i}
  2^{jn(\alpha-1) }
 \int_{\overline {C_j}} |b(y)-b_{\overline {B_j}}|^{r }|f(y)|^rdy
\\
& & +2^{n|\alpha |}2^{n|\alpha -1|}3^{r  }
 \frac{1}{2^{in\alpha }|B(1)|}\sum\limits_{j=-\infty}^{i}
  2^{jn(\alpha-1) }
 \int_{\overline {C_j}} |b_{  {B_i}}-b_{\overline {B_j}}|^{r }|f(y)|^rdy
  \\
&=& 2^{n|\alpha |}2^{n|\alpha -1|}3^{r  }(I_1+I_2+I_3).
\end{eqnarray*}
By (10), we see that
$$
 \overline {C_j}\subset \overline {B_i}=  B(|x|) {\rm ~for~} j\leq i {\rm ~and~}x\in C_i;~~
 {\rm and }~\frac{1}{2^{in }|B(1)|}\leq \frac{1}{|\overline {B_i}|}.
$$
  Let $x\in C_i$, then  we have
 \begin{eqnarray*}
I_1 & \leq &  |b(x)-b_{  {B_i}}|^{r }
  \frac{1}{2^{in }|B(1)|}\int_{  \overline {B_i}}  |f(y)|^rdy
\sum\limits_{j=-\infty}^{i}
  2^{-(i-j)n(\alpha-1) }
  \\
&\leq &  \frac{1}{1-2^{-n(\alpha-1) }}|b(x)-b_{  {B_i}}|^{r }
  \frac{1}{|\overline {B_i}|}\int_{  \overline {B_i}}  |f(y)|^rdy
,
\end{eqnarray*}
noticing that $\alpha-1> 0.$
 For $I_2,$ since $f$ is local bounded, by   Holder's inequality and Lebesque's control convergence
 theorem, we have
 \begin{eqnarray*}
\int_{\overline {C_j}} |b(y)-b_{\overline {B_j}}|^{r }|f(y)|^rdy
& \leq &
\left(\int_{\overline {C_j}} |b(y)-b_{\overline {B_j}}|^{r l'} dy\right)^{1/l'}
\left(\int_{\overline {C_j}}  |f(y)|^{rl}dy\right)^{1/l}
\\
& \leq & |\overline {B_j}|^{1/l'} \|b\|^{ r}_{CMO}\left(\int_{\overline {C_j}}  |f(y)|^{rl}dy\right)^{1/l}
\end{eqnarray*}
\begin{eqnarray}
& ~~~~~~~~~~~~~~~~~~~~~~~~\rightarrow &  \|b\|^{ r}_{CMO}
\int_{\overline {C_j}}  |f(y)|^{r}dy, ~~~({\rm when~}l \rightarrow
1).
 \end{eqnarray}
 Thus, as $I_1$,
 \begin{eqnarray*}
I_2 & \leq &  \|b\|^{ r}_{CMO} \frac{1}{1-2^{-n(\alpha-1) }}   \frac{1}{|\overline {B_i}|}\int_{  \overline {B_i}}  |f(y)|^rdy
.
\end{eqnarray*}
 For $I_3,$ by (10),
\begin{eqnarray*}
 |b_{  {B_i}}-b_{\overline {B_j}}|
&\leq & \frac{1}{|\overline {B_j}|}\int _{\overline
{B_j}}|b(x)-b_{  {B_i}}|dx
\\
& \leq &   \left\{
\begin{array}{cc}
   \frac{1}{|  {B_j}|}\int _{  {B_j}}|b(x)-b_{  {B_i}}|dx,
& {\rm if}~ j\leq i-1,\\
  \frac{1}{|  {B_{j-1}}|}\int _{  {B_j}}|b(x)-b_{  {B_i}}|dx
, & {\rm if}~ j= i,
\end{array}
\right.
\\
&\leq & \frac{2^{n}}{| {B_j}|}\int _{  {B_j}}|b(x)-b_{  {B_i}}|dx
\\
&\leq & \frac{2^{n}}{| {B_j}|}\int _{  {B_j}}|b(x)-b_{  {B_j}}|dx
+2^{n}\sum\limits_{h=j}^{i-1}|b_{  {B_h}}-b_{B_{h+1}}|
\\
 &\leq & 2^{2n}\|b\|_{CMO}(i-j),
\end{eqnarray*}
thus
\begin{eqnarray*}
I_3 & \leq &  \|b\|_{CMO}^{r }
  \frac{1}{2^{in }|B(1)|}\int_{  \overline {B_i}}  |f(y)|^rdy
\sum\limits_{j=-\infty}^{i}
  2^{-(i-j)n(\alpha-1) }(i-j)^{r }
  \\
&\leq &  c\|b\|_{CMO}^{r }
  \frac{1}{|\overline {B_i}|}\int_{  \overline {B_i}}  |f(y)|^rdy
,
\end{eqnarray*}
where $c=  2^{2nr }\sum\limits_{j=-\infty}^{i}
  2^{-(i-j)n(\alpha-1) }(i-j)^{r }.$
    When $x\in C_i$, let
$$  g(x)=\left(
\frac{1}{|\overline {B_i}|}\int_{  \overline {B_i}}
|f(y)|^rdy\right)^{1/r} = \left( \frac{1}{|  {B(|x|)}|}\int_{
{B(|x|)}}  |f(y)|^rdy\right)^{1/r}  ,$$
  combining to
the  estimates of $I_1,I_2,I_3$ above, and noticing that $(|b(x)-b_{
{B_i}}|^{r } +\|b\|^{ r}_{CMO})^{s/r}\leq 2^{s/r}(|b(x)-b_{
{B_i}}|^{s } +\|b\|^{ s}_{CMO})$ for $s>0,r>0$, we see that
\begin{eqnarray}
\int_{C_i}|B(|x|)|^{\gamma -1}
 \left[ h(x)  \right]^{s/r}dx
\leq  c_0 \int_{C_i}|B(|x|)|^{\gamma -1}
 \left[ |b(x)-b_{  {B_i}}|^{s }
+\|b\|^{ s}_{CMO}  \right]^{ }g^s(x)dx
 ,
\end{eqnarray}
where $c_0= [2^{n|\alpha |}2^{n|\alpha -1|}3^{r  }\times4\times 2^{2nr }\sum\limits_{j=-\infty}^{i}
  2^{-(i-j)n(\alpha-1) }(i-j)^{r }]^{s/r}.$
While, since $f$ local bound implies $g$ local bound, as (14),
\begin{eqnarray}
 \int_{C_i}|B(|x|)|^{\gamma -1}
   |b(x)-b_{  {B_i}}|^{s }
 g^s(x)dx
 \leq   \|b\|^{ s}_{CMO}\int_{C_i}|B(|x|)|^{\gamma -1}
  g^s(x)dx.
\end{eqnarray}
Thus, combining to (11), (15) and  (16), we obtain
\begin{eqnarray*}
   M_{s}((M_{r,b  }(f,  \alpha ),\gamma)(R)
&\leq &  c_0^{1/s} \|b\|_{  CMO   }
\left(\frac{1}{|B(R)|^{\gamma}} \sum\limits_{i=-\infty}^{N}
 \int_{C_i}|B(|x|)|^{\gamma -1}
 g^s(x)dx \right)^{1/s}
 \\
& =&    c_0^{1/s} \|b\|_{  CMO   }
M_{s}((M_{r  }(f,  1 ),\gamma)(R).
\end{eqnarray*}
Using (1), we obtain (5).
\par The proof of (6).~~ Replace (9) and (10) by
\begin{eqnarray*}
B_i=\left \{\begin{array}{cc}
B(2^i), & {\rm if}~ i\geq N,\\
B(R), & {\rm if}~ i= N-1,
\end{array}
\right. ,C_i = B_i \setminus B_{i-1}, i=N, N+1, ... , \infty ;
\end{eqnarray*}
 and if $x\in C_i$,
  \begin{eqnarray*}
  {\overline {B_j}}=\left \{
\begin{array}{cc}
B_j, & {\rm if}~ j\geq i ,\\
B(|x|), & {\rm if}~ j= i-1,
\end{array}
\right.
 ,
\overline {C_j} = \overline {B_j} \setminus \overline {B_{j-1}},
j=i, i+1, ..., \infty .
\end{eqnarray*}
Then,
\[
h^*(x)=\left[  (M^*_{r,b  }(f,  \alpha ) (|x|)\right]^r
=\frac{1}{|B(|x|)|^{\alpha}}
\sum\limits_{j=i}^{\infty}
 \int_{\overline {C_j}}|B(|y|)|^{\alpha -1}
(|b(x)-b(y)| |f(y)|)^rdy,
\]
and
\begin{eqnarray*}
  \left[ M^*_{s}((M^*_{r,b  }(f,  \alpha ),\gamma)(R)\right]^s
 = \frac{1}{|B(R)|^{\gamma}} \sum\limits_{i=N}^{\infty}
 \int_{C_i}|B(|x|)|^{\gamma -1}
 \left[ h^*(x)  \right]^{s/r}dx.
\end{eqnarray*}
The rest of the proof of (6) is exactly similar to that of (5).
\par Thus, we finish the proof of Theorem 3.
 \par  The proofs of Theorem 4  by using (5)(6) are exactly similar to that of  Theorem 2 by using (1)(2).

 \par Long Shunchao, Jian Wang,
\par Mathematics Department,
\par Xiangtan University,
\par Xiangtan , 411105,  China,
\par E-mail address: sclong@xtu.edu.cn

\end{document}